\documentclass{amsart}

\usepackage{amsmath,amscd,amssymb,amsthm}
\usepackage{enumerate,mathrsfs}
\usepackage{geometry}
\usepackage{hyperref}
\usepackage{url}
\usepackage{ifthen}
\usepackage[all]{xy} \xyoption{2cell} \UseAllTwocells

\numberwithin{equation}{section}

\theoremstyle{plain}
\newtheorem{thm_}[equation]{Theorem}
\newtheorem{lemma_}[equation]{Lemma}
\newtheorem{prop_}[equation]{Proposition}
\newtheorem{cor_}[equation]{Corollary}
\newtheorem{eg_}[equation]{Example}
\newtheorem{con_}[equation]{Conjecture}
\newtheorem*{cons_}{Conjecture}

\theoremstyle{definition}
\newtheorem{thmu_}[equation]{Theorem}
\newtheorem*{thmus_}{Theorem}
\newtheorem{propu_}[equation]{Proposition}
\newtheorem*{propus_}{Proposition}
\newtheorem{coru_}[equation]{Corollary}
\newtheorem*{corus_}{Corollary}
\newtheorem{lemu_}[equation]{Lemma}
\newtheorem*{lemus_}{Lemma}
\newtheorem{egu_}[equation]{Example}
\newtheorem*{egus_}{Example}
\newtheorem{def_}[equation]{Definition}
\newtheorem*{defs_}{Definition}
\newtheorem{rk_}[equation]{Remark}
\newtheorem*{rks_}{Remark}
\newtheorem{ex_}[equation]{Remark}

\newcommand{\thm}[1]{\begin{thm_}#1\end{thm_}}
\newcommand{\thmu}[1]{\begin{thmu_}#1\end{thmu_}}

\newcommand{\lemm}[1]{\begin{lemma_}#1\end{lemma_}}

\newcommand{\prop}[1]{\begin{prop_}#1\end{prop_}}

\newcommand{\defi}[1]{\begin{def_}#1\end{def_}}

\newcommand{\rk}[1]{\begin{rk_}#1\end{rk_}}

\newcommand{\cor}[1]{\begin{cor_}#1\end{cor_}}

\newcommand{\pf}[1]{\begin{proof}#1\end{proof}}

\DeclareMathOperator{\Hom}{Hom}

\DeclareMathOperator{\Spec}{Spec}

\DeclareMathOperator{\im}{im}

\DeclareMathOperator{\Br}{Br}

\newcommand{\bA}{\mathbb A}%

\newcommand{\bfA}{\mathbf A}%
\newcommand{\bfG}{\mathbf G}%

\newcommand{\cS}{\mathcal S}

\newcommand{\cV}{\mathcal V}
\newcommand{\cW}{\mathcal W}
\newcommand{\cX}{\mathcal X}
\newcommand{\cY}{\mathcal Y}
\newcommand{\cZ}{\mathcal Z}%

\newcommand{\s}{\sigma}

\newcommand{\tm}{\times}%

\newcommand{\ra}{\rightarrow}

\newcommand{\xra}{\xrightarrow}

\newcommand{\mpt}{\mapsto}

\newcommand{\os}[2]{\overset{#1}{#2}}

\newcommand{\is}[2]{\xymatrix@-4mm{#1 \ar[r]^-{\sim} & #2 }}
\newcommand{\mis}[2]{\xymatrix@-2mm{#1 \ar[r]^-{\sim} & #2 }}

\newcommand{\dra}[4]{\xymatrix@-4mm{#1 \ar@<.5ex>[r]^-{#3} \ar@<-.5ex>[r]_-{#4}& #2 }}
\newcommand{\era}[5]{\xymatrix@-4mm{#1 \ar[r] &#2 \ar@<.5ex>[r]^-{#4} \ar@<-.5ex>[r]_-{#5}& #3 }}

\newcommand{\tu}[1]{\text{\upshape #1}}

\newcommand{\SSet}{{Set}}

\newcommand{\SSch}{{Sch}}

\newcommand{\op}{\tu{op}}

\newcommand{\et}{\tu{\'et}}

\newcommand{\fppf}{\tu{fppf}}

\DeclareFontFamily{U}{wncy}{}
\DeclareFontShape{U}{wncy}{m}{n}{%
   <5>wncyr5%
   <6>wncyr6%
   <7>wncyr7%
   <8>wncyr8%
   <9>wncyr9%
   <10>wncyr10%
   <11>wncyr10%
   <12>wncyr6%
   <14>wncyr7%
   <17>wncyr8%
   <20>wncyr10%
   <25>wncyr10}{}
\DeclareMathAlphabet{\cyrille}{U}{wncy}{m}{n}
\def\Sha{\cyrille X}

\newcommand{\Tors}{\mathsf{Tors}}

\newcommand{\eq}[1]{\begin{equation}#1\end{equation}}
\newcommand{\eqn}[1]{\begin{equation*}#1\end{equation*}}

\newcommand{\gan}[1]{\begin{gather*}#1\end{gather*}}

\newcommand{\aln}[1]{\begin{align*}#1\end{align*}}

\newcommand{\aci}[1]{\ar@{^(->}[#1]|-{/}}
\newcommand{\coaci}[1]{\ar@{_(->}[#1]|-{/}}
\newcommand{\aoi}[1]{\ar@{^(->}[#1]|-{\circ}}
\newcommand{\coaoi}[1]{\ar@{_(->}[#1]|-{\circ}}
\makeatletter
\def\citet@url@sp{https://stacks.math.columbia.edu/}
\def\citet@bib@sp{stacks-project}
\def\citet@url@kd{https://kerodon.net/}
\def\citet@bib@kd{kerodon}
\newcommand{\citet@tag}[2]{\href{#2tag/#1}{#1}}
\newcommand{\citet@taglist}[2]{%
 \def\@citet@e{}%
 \def\@citet@tag@n{0}
 \@for\@citet@tag:=#1\do{%
  \edef\@citet@tag@n{\the\numexpr\@citet@tag@n + 1}%
 }%
 \def\@citet@tags{%
  \def\@citet@tag@i{0}%
  \@for\@citet@tag:=#1\do{%
   \edef\@citet@tag@i{\the\numexpr\@citet@tag@i + 1}%
   \ifthenelse{\@citet@tag@i > 1}{
    \ifthenelse{\@citet@tag@i = \@citet@tag@n}{
     \citet@seplast%
    }{%
     \citet@sep%
    }%
   }{}%
   \citet@entry{\@citet@tag}{#2}
  }%
 }%
 \ifthenelse{\@citet@tag@n > 1}{
  \def\@citet@Tag{Tags}%
 }{%
  \def\@citet@Tag{Tag}%
 }%
 \@citet@Tag~\@citet@tags
}
\newcommand{\citet@sep}{, }
\newcommand{\citet@seplast}{ and }
\newcommand{\citet@entry}[2]{\citet@tag{#1}{#2}}
\let\@old@cite\cite
\renewcommand{\cite}[2][]{%
 \def\@citet@detail{\citet@taglist{#1}{\@citet@url}}%
 \ifthenelse{\equal{#2}{sp}}{%
  \def\@citet@url{\citet@url@sp}%
  \def\@citet@bib{\citet@bib@sp}%
 }{\ifthenelse{\equal{#2}{kd}}{%
  \def\@citet@url{\citet@url@kd}%
  \def\@citet@bib{\citet@bib@kd}%
 }{
  \def\@citet@detail{#1}%
  \def\@citet@bib{#2}%
 }}%
 \ifthenelse{\equal{#1}{}}{%
  \@old@cite{\@citet@bib}%
 }{%
  \@old@cite[\@citet@detail]{\@citet@bib}%
 }%
}
\makeatother

\newcommand{\etale}{{\'etale}}

\newcommand{\BM}{{Brauer-Manin}}

\newcommand{\adelic}{{ad\`elic}}

\newcommand{\desc}{\tu{desc}}
\newcommand{\ddesc}{{\desc, \desc}}

\newcommand{\etBr}{{\et, {\Br}}}

\newcommand{\hdesc}[2]{{#1\tu{-}\desc{\ifthenelse{\equal{#2}{}}{}{_#2}}}}

\newcommand{\XA}{X(\bfA_k)}

\newcommand{\Xk}{X(k)}

\newcommand{\cXAk}{\cX(\bfA_k)}

\newcommand{\fin}{\tu{fin}}
\newcommand{\fdesc}{{\fin, \desc}}
\newcommand{\ob}{\tu{ob}}

\begin{document}
\title[Finite descent for stacks]{Iterated descent obstructions for algebraic stacks}

\author[H. Wu]{Han Wu}
\address{Hubei Key Laboratory of Applied Mathematics,
	Faculty of Mathematics and Statistics,
	Hubei University,
	No. 368, Friendship Avenue, Wuchang District, Wuhan,
	Hubei, 430062, P.R.China.}
\email{wuhan90@mail.ustc.edu.cn} 

\author[C. Lv]{Chang Lv}
\address{Key Laboratory of Cyberspace Security Defense\\
		Institute of Information Engineering\\
		Chinese Academy of Sciences\\
		Beijing 100093, P.R. China}
\email{lvchang@amss.ac.cn}

\date{\today}
\thanks{C. Lv is partially supported by
		National Natural Science Foundation of China NSFC Grant No. 11701552. H. Wu is
		partially supported by NSFC Grant No. 12071448.}

\begin{abstract}
		Base on a conjecture, we prove that for any smooth separated stack of finite type over a number field, its descent obstruction equals its iterated descent obstruction.  As a consequence, we show that
		for any algebraic stack over a number field
    that has a finite \'etale covering of a
		smooth geometrically integral variety, its
		descent obstruction equals its iterated descent obstruction.
\end{abstract}
\maketitle

\section{Introduction} \label{intro}
The existence of rational points is a fundamental arithmetic property of
 varieties over number fields.
The most common principle is the local-global principle.
When it fails, one considers various obstructions to it, such as the
 {\BM} obstruction and the descent obstruction.
We  start by reviewing the classical definitions and results.

Let $X$  be a variety over a number field $k$, we write $X(-)$ for the functor
 $\Hom_k(-, X)$. The \emph{rational points} of $X$ is the set $\Xk$, which is
 contained in its \emph{ad\`elic points} set $\XA$.
We say that the \emph{local-global principle} holds if $\XA\neq\emptyset$
 implies $\Xk\neq\emptyset$.

Let $F\colon (\SSch/k)^\op\ra \SSet$ be a contravariant
 functor from the category of $k$-schemes
 to the category of sets.
For a $k$-scheme $T$ and $A\in F(X)$, the \emph{evaluation} of
 a \emph{$T$-point} $T\os{x}{\ra}X \in X(T)$ at
 $A$ is defined to be  the image  of  $A$ under the pull-back map
 $F(x)\colon F(X)\ra F(T)$  induced by $x$, denoted by $A(x)$.
We have an obvious commutative diagram
\eq{ \label{eq_F_ob}
\xymatrix{
\Xk\ar[r]\ar[d]^-{A(-)} &\XA\ar[d]^-{A(-)} \\
F(k)\ar[r] &F(\bfA_k)
}}
Define the \emph{obstruction $($set$)$ given by $A$} to be
  $\XA^A=\{x\in \XA\mid A(x)\in \im(F(k)\ra F(\bfA_k)\}$. Then we
 have $ \Xk\subseteq \XA^A\subseteq \XA$, giving
 a constraint on the locus of  rational points in {\adelic} points.
Putting
\eqn{
\XA^F\colon=\XA^{F(X)}=\bigcap_{A\in F(X)} \XA^A,
}
 called the \emph{$F$-set} or \emph{$F$-obstruction},
 which also yields inclusions
\eqn{
\Xk\subseteq \XA^F\subseteq\XA.
}
See Poonen \cite[8.1.1]{poonen17rational}.

\defi{
We say that the $F$-obstruction to the local-global principle
 is the \emph{only one} if
 $\XA^F\neq\emptyset$ implies $\Xk\neq\emptyset$.

We say that \emph{there is} a
 $F$-obstruction to the local-global principle if $\XA\neq\emptyset$ but
 $\XA^F=\emptyset$ (a \emph{priori} $\Xk=\emptyset$).
}

We are mainly interested in \emph{cohomological obstructions}, namely, ones
 where $F$ is taken to be a cohomological functor.

Let $F=\Br\colon=H^2_{\text{\'et}}(-, \bfG_m)$ be
 the cohomological Brauer-Grothendieck group (c.f. Grothendieck
 \cite{grothendieck95brauer}).
Thus we obtain $X(\bA_k)^{\Br}$, the {\BM} obstruction (see, e.g. Skorobogatov
 \cite{torsor}).

Another cohomological functor is $F=\check H_\fppf^1(-, G)$, the obstruction given by
 the first \v Cech cohomology, where $G$ is a linear $k$-group. In general,
 $\check H_\fppf^1(X, G)$ is a pointed set, which is isomorphic to
 $H_\fppf^1(X, G)$ if $G$ is commutative, and further  to $H_\et^1(X, G)$ since
 $G$ is smooth over the number field $k$.
The classical \emph{descent obstruction} is given by
\eqn{ \label{eq_cdesc}
\XA^\desc\colon=
 \bigcap_{\text{all linear $k$-group $G$}} \XA^{\check H_\fppf^1(X, G)}.
}
The descent theory was established by Colliot-Th{\'e}l{\`e}ne  and Sansuc
 \cite{cs87descente-ii} for tori  and Skorobogatov
 \cite{skorobogatov99beyond} for groups of multiplicative type.
Harari \cite{harari02groupes}, Harari
 and Skorobogatov \cite{hs02non-abelian,
 hs05non-abelian} studied the descent obstruction for general algebraic groups
 and compared it with the {\BM} obstruction.
One of the results is the well-known inclusion
 $\XA^\desc\subseteq \XA^{\Br}$ for regular, quasi-projective $k$-variety $X$
 (see, e.g., \cite[Prop. 8.5.3]{poonen17rational}).

The classical \emph{descent by torsors} says that
\eqn{
\XA^f\colon=\bigcup_{\s\in  H^1(k, G)} f_\s(Y_\s(\bfA_k))
}
 where $f\colon Y\os{G}{\ra} X$ is an $X$-torsor under the linear $k$-group $G$,
 representing the class $[Y]\in \check H_\fppf^1(X, G)$,
 and $f_\s\colon Y_\s\os{G_\s}{\ra} X$ is the twisted torsor.
Then we may define various of composite (or iterated) obstruction sets
 $\XA^{\desc, \ob}$ between $\Xk$ and $\XA^\desc$, such as
 the \emph{iterated descent obstruction}
\eqn{
\XA^\ddesc=\colon \bigcap_{\substack{
   \text{all linear $G$} \\ \text{all torsor } f\colon Y\os{G}{\ra} X}}
 \bigcup_{\s\in  H^1(k, G)} f_\s(Y_\s(\bfA_k)^\desc)
}
 and the \emph{{\etale} Brauer obstruction}
\eqn{
\XA^\etBr\colon=
 \bigcap_{\substack{
   \text{all finite $G$} \\ \text{all torsor } f\colon Y\os{G}{\ra} X}}
 \bigcup_{\s\in  H^1(k, G)} f_\s(Y_\s(\bfA_k)^{\Br}).
}
By works of Stoll, Skorobogatov, Demarche, Poonen, Xu and Cao,
 it is known that for $X$ being a
 smooth, quasi-projective, geometrically integral
 $k$-variety,  $\XA^\desc=\XA^\etBr$ \cite{stoll07finite,
 skorobogatov09descent, demarche09obstruction, poonen10insufficiency,
 cdx19comparing} and
\eq{ \label{eq_sch_dd=d}
\XA^\ddesc=\XA^\desc
}
 \cite{cao20sous} and no smaller
 obstruction set than the descent set is discovered.

In this paper, all stacks are smooth separated stack of finite type over a number field. We consider  the
 descent obstruction $\cXAk^{\desc}$ for an algebraic $k$-stack $\cX$. Our main result base on a conjecture: similar to the existence of Galois closure of the composite of two Galois extensions, we conjecture that there exists a torsor covering of the composite of two torsors, see Conjecture \ref{main conj} for more details.
Our main result are the following.
\thm{ [Theorem \ref{main theorem dd=d conj}]
	Let $\cX$ be an algebraic stack over a number field $k$. Assume Conjecture \ref{main conj} holds.
	Then 
	\eqn{
		\cXAk^\desc = \cXAk^\ddesc.
}}
Without any assumption, we have the following theorem.

\thmu{[Theorem \ref{thm_dd=d}]
			Let $\cX$ be an algebraic stack over a number field $k$. Assume there exists a \'etale finite covering $f\colon Y\to X$, where $Y$ is a quasi-projective smooth
		geometrically integral $k$-variety.
		Then we have
		\eqn{
			\cXAk^\desc=\cXAk^\ddesc.
}}
This extend \eqref{eq_sch_dd=d} from varieties to algebraic stacks having finite \'etale covering by varieties. The main theorem are proven in Section \ref{main}.

\section{Composite of torsors} \label{cover}
We consider the composite of two torsors. It is crucial to our paper.

Let $f\colon \cY\xra{F} \cX$  and
Let $g\colon \cZ\xra{G} \cY$  be two torsors of algebraic $k$-stacks
 (see \cite[Sec. 2.3]{lh23stackbm}), where $F$
	and $G$ are linear $k$-groups.
If $F$ and $G$ are finite groups, then by the existence of Galois covering of $f\circ g\colon \cZ\to \cX$, there exist some finite $k$-group $H$
and an $H$-torsor $h\colon \cV\xra{H} \cX$ such that the following
diagram is commutative and all arrows are surjective torsor morphisms
\[\xymatrix{
	\cV\ar[d]\ar[dr]\ar^{H}[drr] & & \\
	\cZ\ar_{G}[r] &\cY\ar_{F}[r] &\cX.}  \]	     	
     If $F$ and $G$ are not finite groups, the  existence of some linear $k$-group $H$ having above properties seems to be difficult. But we strongly believe this is true, and we post it as a conjecture.
     \begin{con_}\label{main conj}
     	Given two linear $k$-groups  $F$ and $G$, let  $f\colon \cY\xra{F} \cX$  and
     	 $g\colon \cZ\xra{G} \cY$  be two torsors of algebraic $k$-stacks, then there exist some linear $k$-group $H$
     	 and an $H$-torsor $h\colon \cV\xra{H} \cX$ such that the following
     	 diagram is commutative and all arrows are surjective torsor morphisms
     	 \[\xymatrix{
     	 	\cV\ar[d]\ar[dr]\ar^{H}[drr] & & \\
     	 	\cZ\ar_{G}[r] &\cY\ar_{F}[r] &\cX.}  \]	     	
     \end{con_}
 Though in the case that $\cX,\cY,\cZ$ are smooth, geometrically integral varieties, this conjecture is unknown. In this paper,
 we will give an explicit construction of $\cV$ in a special case.
Since this construction rely on the existence of Weil restriction.
Firstly, we consider Weil restriction to algebraic stacks.
\defi{
		Let $\cZ\ra\cY\xra{f}\cX$ be morphisms of algebraic $k$-stacks. If the functor
\aln{
			\SSch/\cX&\ra \SSet \\
			\cS&\mpt \Hom_\cY(\cS\tm_\cX \cY, \cZ)
}
		is represented by a stack,
    then we call it the \emph{Weil restriction} and denote  it by $R_{\cY/\cX}\cZ$.
}

For the morphism $f\colon \cY\to\cX$ is finite, the Weil restriction
 $\cW=R_{\cY/\cX}\cZ\ra \cX$ exists as algebraic
 $k$-stack (c.f. \cite[Cor. 9.2 (ii)]{hr19algebraicity} in which
  $\cW$ is denoted by  $f_*\cZ$ and called the pushforward along $f$).
We assume that $\cY\to \cX$ is an $F$-torsor and $F$ is a finite $k$-group.
Then we have the following
\prop{ \label{prop_tor_cov_fin}
Let $f\colon \cY\xra{F} \cX$  and
Let $g\colon \cZ\xra{G} \cY$  be two torsors of algebraic
 $k$-stacks, where $F$ is a finite $k$-group
		and $G$ is an arbitrary linear $k$-group.
Then Conjecture \ref{main conj} holds, i.e. there exist some linear $k$-group $H$ and a torsor
 $h\colon \cV\xra{H} \cX$ such that the following diagram is commutative
 and all arrows are surjective torsor morphisms.
		$$
		\xymatrix{
			\cV\ar[d]\ar[dr]\ar^{H}[drr] & & \\
			\cZ\ar_{G}[r] &\cY\ar_{F}[r] &\cX.
		}$$
}

\pf{
The proof is along the same lines as the proof of
 \cite[Prop. 2.3]{skorobogatov09descent}, where the same
 statement is shown for varieties.
Hence, we sketch it by focusing on their differences.
Since $F$ is finite, the Weil restriction
		$\cW=R_{\cY/\cX}\cZ\ra \cX$ exists as algebraic
		$k$-stack,
    and it is a torsor
		under the group $R_{\cY/\cX}(G_\cY)$.
Note that $R_{F/k}(G_F)$ is a linear $k$-group and
		\gan{
			R_{\cY/\cX}(G_\cY)\tm_\cX\cY \xra{\sim} R_{F/k}(G_F)\tm_k \cY \text{ and } \\
			R_{\cY/\cX}(G_\cY)\tm_\cX F_\cX \xra{\sim} R_{F/k}(G_F)\rtimes F.
		}
We see that $\cV\colon=\cW\tm_\cX \cY$ is in both $\Tors(\cY, R_{F/k}(G_F))$ and
		$\Tors(\cX, R_{F/k}(G_F)\rtimes F)$.
Finally, the map $\cV\ra \cZ$ is obtained by the unit $\colon\Spec k\ra F$.
}

\section{Iterated descent obstructions} \label{descent}
Iterated descent obstruction is already defined for schemes, cf. \cite[Chapter 8]{poonen17rational}. 
Since for a (resp. Commutative) $k$ group $G$,
$\check H_\fppf^1(-, G)$ (resp. $H_\et^2(-, G)$) is  stable
(see \cite[2.17 Def.]{lv2desc}),
the following definitions
also make sense by the diagram for $\cX$ similar to \eqref{eq_F_ob}.
\defi{ \label{defi_desc_conn_2desc}
	The \emph{descent obstruction} is
	\eqn{
		\cXAk^\desc=\bigcap_{\text{all linear $k$-groups  $G$}}
		\cXAk^{\check H_\fppf^1(-, G)}.
	}
	Also, descent along a torsor is also correct for algebraic stacks
	\cite[3.20 Prop.]{lv2desc},
	we may define various of composite (or iterated) obstruction sets
	\eqn{
		\cXAk^\fdesc=
		\bigcap_{\substack{\text{all finite $F$ } \\
				\text{all torsor } f\colon \cY\os{F}{\ra} \cX}}
		\bigcup_{\s\in  H^1(k, G)} f_\s(\cY_\s(\bfA_k)^\desc),
	}
	\eqn{
		\cXAk^\ddesc=
		\bigcap_{\substack{\text{all linear $G$ } \\
				\text{all torsor } f\colon \cY\os{G}{\ra} \cX}}
		\bigcup_{\s\in  H^1(k, G)} f_\s(\cY_\s(\bfA_k)^\desc).
	}
}

The following lemma will be used in the proof of Theorem \ref{prop_fd_d}. We prove it first.
\lemm{ \label{lemma:the closure of F and G plus conj}
Assume Conjecture \ref{main conj} holds. Let $F$  and $G$ be a linear algebraic group.
Let  $f\colon \cY\ra \cX$ be an $F$-torsor.
Let  $P\in f(\cY(\bfA_k)).$ If $P\notin f(\cY(\bfA_k)^{H_\fppf^1(\cY, G)})$,
  then there exist linear algebraic groups $G',\tilde{G}$ and a $\cY$-torsor
  $g\colon \cZ\to \cY$ under the group $G'$ such that:
\begin{enumerate}
			\item\label{lemma cond1 conj} $G$ is a quotient group of $G'$,
			\item\label{lemma cond2 conj} $\cZ$ is an $\cX$-torsor under the group $\tilde{G}$,
			\item\label{lemma cond3 conj} $1\to G'\to \tilde{G}\to F\to 1$ is an exact sequence of groups, and
			$$	\xymatrix{
				\cZ\ar_{G'}[d]\ar^{\tilde{G}}[dr] &  \\
				\cY\ar_{F}[r] &\cX.
			}$$
			the actions of $F$, $G'$ and $\tilde{G}$ are compatible,
			\item  $$P\notin f(\bigcup_{\s\in H^1(k,G')}g_\s(\cZ_\s(\bfA_k))).$$
      In particular, $P\notin f(\cY(\bfA_k)^{H_\fppf^1(\cY, G')}).$
\end{enumerate}

}
\pf{
By the assumption that $P\notin f(\cY(\bfA_k)^{H_\fppf^1(\cY, G)})$,
 there exists a $\cY$-torsor $g'\colon \cZ'\to \cY$ under the group $G$
 such that
\eqn{
P\notin f(\bigcup\limits_{\s'\in H^1(k,G)}g'_{\s'}(\cZ'_{\s'}(\bfA_k))).
}
By assumption that Conjecture \ref{main conj} holds, there exist a linear $k$-group
 $\tilde{G}$ and a torsor $ \cZ\xra{\tilde{G}} \cX$ such that the
 following diagram is commutative and all arrows are surjective torsor morphisms.
		$$
		\xymatrix{
			\cZ\ar[d]\ar[dr]\ar^{\tilde{G}}[drr] & & \\
			\cZ'\ar_{G}[r] &\cY\ar_{F}[r] &\cX.
		}$$
Let $G'$ is the kernel of the surjective morphism $\tilde{G}\to F$,
 then the conditions \eqref{lemma cond1 conj}, \eqref{lemma cond2 conj} and
 \eqref{lemma cond3 conj} are satisfied.

We use contradiction to prove that
 $P\notin f(\bigcup\limits_{\s\in H^1(k,G')}g_\s(\cZ_\s(\bfA_k))).$
Otherwise, there exists some $\s_0\in H^1(k,G')$ such that
 $P\in f(g_{\s_0}(\cZ_{\s_0}(\bfA_k)))$.
Let $\s'_0\in  H^1(k,G')$ be the image of $\s_0$ under the natural map
 $H^1(k,G')\to H^1(k,G)$.
We have the following diagram:
		$$
		\xymatrix{
			\cZ_{\s_0}\ar[d]\ar^{G'_{\s_0}}[dr] &  \\
			\cZ'_{\s'_0}\ar_{G_{\s'_0}}[r] &\cY,
		}$$
	here ${G'_{\s_0}}$ and ${G_{\s'_0}}$ are inner of $G'$ and $G$ given
  by $\s_0$ and $\s'_0$ respectively.
Then
\eqn{
P\in f(g_{\s_0}(\cZ_{\s_0}(\bfA_k)))\subset
		f(g'_{\s'_0}(\cZ'_{\s'_0}(\bfA_k))),
}
 which is a contradiction.
Hence  $P\notin f(\bigcup\limits_{\s\in H^1(k,G')}g_\s(\cZ_\s(\bfA_k)))$.
}

\rk{\label{remark: the closure of F and G}
	Assume Conjecture \ref{main conj} holds.
If we have a  $\cY$-torsor $g'\colon \cZ'\to \cY$ under the group
 $G$ such that
\eqn{
P\notin f(\bigcup\limits_{\s\in H^1(k,G)}g'_\s(\cZ'_\s(\bfA_k))),
}
 then we replace $G$ by $G'$ and $\cY$-torsor $\cZ'$ by
 $g\colon \cZ\to \cY$ under the group $G'$ such that
 $P\notin f(\bigcup\limits_{\s\in H^1(k,G)}g_\s(\cZ_\s(\bfA_k)))$
 and the composite of $\cZ\to \cY\to\cX$ is a
		$\cX$-torsor.
		In particular, if $F$ is a finite group, then by Proposition \ref{prop_tor_cov_fin}, Conjecture \ref{main conj} holds, and we can do the same thing.
}

\section{Main result} \label{main} 

With the preparation, we are able to prove our main theorem.

\thm{ \label{main theorem dd=d conj}
		Let $\cX$ be an algebraic stack over a number field $k$. Assume Conjecture \ref{main conj} holds.
		Then 
		\eqn{
			\cXAk^\desc = \cXAk^\ddesc.
}}
\pf{
By the functoriality of descent obstructions,
 we have
\eqn{\cXAk^\desc \supset \cXAk^\ddesc.
}
Next, we prove the other side:
\eqn{
\cXAk^\desc \subset \cXAk^\ddesc.
}

We take a point $P\in \cXAk^\desc$.
For any linear algebraic group $F$ and any $\cX$-torsor $f\colon \cY\to \cX$ under the group $F$, we need to prove that $P\in \bigcup\limits_{\s\in H^1(k,F)} f_\s(\cY_\s(\bfA_k)^\desc).$

By the definition of descent obstruction,
 there exists an element $\s'\in  H^1(k,F)$ such that
 $P\in f_{\s'}(\cY_{\s'}(\bfA_k))$.
Then over $\bfA_k$, $f_{\s'}^{-1}(P)\cong F_{\s'}$.
Here $F_{\s'}$ is the inner form of $F$ given by ${\s'}$.
Since the set
 $$\Sha(F/K)\colon=\ker(H^1(K,F)\to \prod_{v\in \Omega_K} H^1(K_v,F))$$
 is a finite set \cite[Thm.  5.12.29]{poonen17rational},
 and the equality
\eqn{
\sharp \Sha(F/K)=\sharp \{\s\in  H^1(k,F)| P\in f_{\s}(\cY_{\s}(\bfA_k))\}
}
 holds,
		the set $\{\s\in  H^1(k,F)| P\in f_{\s}(\cY_{\s}(\bfA_k))\}$
    having the same number as the set $\Sha(F/K)$  is also a finite set.
We assume
		that this set equals $\{\s_1,\cdots, \s_m\}$.
We use contradiction to prove the other side.
We assume that $$P \notin \bigcup_{\s\in H^1(k,F)} f_\s(\cY_\s(\bfA_k)^\desc).$$

Then for any integer $i\in \{1,\cdots,m\}$,
 we have $P\notin f_{\s_i}(\cY_{\s_i}(\bfA_k)^\desc)$.
Hence there exists a linear algebraic $k$-group $G_i$
 such that $P\notin f_{\s_i}(\cY_{\s_i}(\bfA_k)^{H_\fppf^1(\cY_{\s_i}, G_i)})$.
By Lemma \ref{lemma:the closure of F and G plus conj},
 we may replace $G_i$ by some linear algebraic group such that
 there exists a $\cY_{\s_i}$-torsor
 $g_i\colon \cZ_i\to \cY_{\s_i}$ under the group $G_i$,
 also an $\cX$-torsor $f_{\s_i} \circ g_i\colon\cZ_i\to \cX$
 under some linear algebraic group $\tilde{G_i}$ (denoted by $\tilde{f_i}$),
 such that
 $$P\notin f_{\s_i}(\bigcup\limits_{\tau\in H^1(k,G_i)} g_{i\tau}(\cZ_{i\tau}(\bfA_k))).$$
We have an exact sequence $1\to G_i\to \tilde{G_i}\to F_{\s_i}\to 1$.
Since $P\in \cXAk^\desc$,
 the set
\eqn{
\{\tilde{\s}\in  H^1(k,\tilde{G_1}) | P\in \tilde{f_1}_{\tilde{\s}}
 (\cZ_{1\tilde{\s}}(\bfA_k)))\}
}
 is nonempty.
Consider its image in $H^1(k,F_{\s_1}),$ denoted by $S_1$.
Since the unit element is not in its image,  we have $\sharp S_1<m$.
Take an $\tilde{\s}\in  H^1(k,\tilde{G_1})$ such that
 $P\in \tilde{f_1}_{\tilde{\s}}(\cZ_{1\tilde{\s}}(\bfA_k)))$,
 and assume that its image in $S_1$ is $\s$.
Then we have a surjective group homomorphism
 $(\tilde{G_1})_{\tilde{\s}}\to (F_{\s_1})_\s$,
 denoted the kernel by ${G_1}_{\tilde{\s}}$.
Since
\eqn{
P\in \tilde{f_1}_{\tilde{\s}}(\cZ_{1\tilde{\s}}(\bfA_k)))\subset
 (f_{\s_1})_\s((\cY_{\s_1})_\s(\bfA_k)),
}
 without loss of generality, we assume $(F_{\s_1})_\s=F_{\s_2}$.
Consider the $(\cY_{\s_1})_\s=\cY_{\s_2}$-torsor
 $\cZ_{1\tilde{\s}}\times_{\cY_{\s_2}}\cZ_2$ under the group
 ${G_1}_{\tilde{\s}}\times G_2$.
Then it is a
		 $\tilde{G_1}_{\tilde{\s}}\times_{F_{\s_2}} \tilde{G_2}$-torsor
     over $\cX$.
We have an exact sequence $$1\to {G_1}_{\tilde{\s}}\times G_2\to
		  \tilde{G_1}_{\tilde{\s}}\times_{F_{\s_2}} \tilde{G_2}\to F_{\s_2}\to 1.$$
Replace $G_2$ and $\cZ_2$ by ${G_1}_{\tilde{\s}}\times G_2$ and
		   $\cZ_{1\tilde{\s}}\times_{\cY_{\s_2}}\cZ_2$ respectively.
Since $P\in \cXAk^\desc$,
 the set
\eqn{
\{\tilde{\s}\in  H^1(k,\tilde{G_2}) | P\in
 \tilde{f_2}_{\tilde{\s}}(\cZ_{2\tilde{\s}}(\bfA_k)))\}
}
 is nonempty.
Consider its image in $H^1(k,F_{\s_2}),$ denoted by $S_2$.
Then $\sharp S_2<\sharp S_1<m$.
By induction, we can get a strictly decreasing sequence. This is impossible.
Hence
\eqn{
\cXAk^\desc = \bigcup\limits_{\s\in H^1(k,F)} f_\s(\cY_\s(\bfA_k)^\desc).
} 	\eqn{
\cXAk^\desc = \cXAk^\ddesc.
}
}

\rk{
Note that unlike \cite{skorobogatov09descent,cdx19comparing},
 our proof does not make use of the compactness argument for {\adelic} topology,
 which is easier to generalize to algebraic stacks.
}


\prop{ \label{prop_fd_d}
	Let $k$ be a number field, and $\cX$ an algebraic $k$-stack.
	Let $F$ be a finite $k$-group and  $f\colon \cY\ra \cX$ be an $F$-torsor.
	Then we
	\eqn{
		\cXAk^\desc = \bigcup_{\s\in H^1(k,F)} f_\s(\cY_\s(\bfA_k)^\desc).
}}

\pf{By the last part of Remark \ref{remark: the closure of F and G}, and the argument as in the proof of Theorem \ref{main theorem dd=d conj}, this proposition follows.
}

Without any assumption, we have the following theorem.
\thmu{\label{thm_dd=d}
	Let $\cX$ be an algebraic stack over a number field $k$. Assume there exists a finite \'etale covering $ X\to \cX$, where $Y$ is a quasi-projective smooth
	geometrically integral $k$-variety.
	Then we have
	\eqn{
		\cXAk^\desc=\cXAk^\ddesc.
}}

\pf{By applying the same argument as in the proof of \cite[Proposition 5.3.9]{Sz08} to stacks, we
	let $f\colon Y\ra \cX$ be the \'etale Galois covering  of  $ X\to \cX$ with the Galois group $F$. By Proposition
	\ref{prop_fd_d}, we have
	\eqn{
		\cXAk^\desc = \bigcup_{\s\in H^1(k,F)} f_\s(Y_\s(\bfA_k)^\desc).
	}
	By \cite[Thm. 1.2]{cao20sous},  $Y_\s(\bfA_k)^\desc=Y_\s(\bfA_k)^\ddesc$.
	Thus  by functoriality
	\eqn{
		\cXAk^\desc = \bigcup_{\s\in H^1(k,F)} f_\s(Y_\s(\bfA_k)^\ddesc) \subseteq
		\cXAk^\ddesc.
	}
	Then other inclusion is clear and the proof is complete.
}
 
 With this theorem, we have the following corollary.
\cor{ \label{prop_fd=d}
	Let $k$ be a number field.
	Let $\cX$ be an algebraic $k$-stack.
	Then we have
	\eqn{
		\cXAk^\desc=\cXAk^\fdesc.
}}
\pf{ It is clear $\cXAk^\desc\subset \cXAk^\fdesc\subset \cXAk^\ddesc.$ 
By Proposition \ref{prop_fd_d}, this corollary follows.
}

\section*{Acknowledgment}
The authors would like to thank Yang Cao and Lei Zhang for helpful discussions.
The work was partially done during the first author's  visits to
the Morningside Center of Mathematics, Chinese
Academy of Sciences. They  thank the center for its hospitality.

\providecommand{\bysame}{\leavevmode\hbox to3em{\hrulefill}\thinspace}
\providecommand{\MR}{\relax\ifhmode\unskip\space\fi MR }
\providecommand{\MRhref}[2]{%
	\href{http://www.ams.org/mathscinet-getitem?mr=#1}{#2}
}
\providecommand{\href}[2]{#2}

\end{document}